\documentclass[utf8]{elsarticle}
\usepackage{lineno}
\usepackage{natbib}
\usepackage{xcolor}
\usepackage{comment}
\usepackage{amsmath}
\usepackage{natbib}
\usepackage{graphicx}
\usepackage[section]{placeins}
\usepackage[colorlinks=true,  linkcolor=blue, citecolor=blue, bookmarksopen=true]{hyperref}
\journal{Elsevier}
\date{}
\begin{document}
\begin{frontmatter}
\title{A multi-period game-theoretic approach to market fairness in oligopolies}

\author[1]{Asimina Marousi}
\affiliation[1]{Department of Chemical Engineering, Sargent Centre for Process Systems Engineering, University College London, Torrington Place, London WC1E 7JE, United Kingdom}
\author[2]{Karthik Thyagarajan}
\author[2]{Jose M. Pinto}
\author[1]{Lazaros G. Papageorgiou}
\author[1]{Vassilis M. Charitopoulos\corref{cor1}}
\ead{v.charitopoulos@ucl.ac.uk}
\affiliation[2]{Linde Digital Americas, 10 Riverview Drive, Danbury CT 06810, United States}
\cortext[cor1]{corresponding author}

\begin{abstract}
Contemporary process industries are constantly confronted with volatile market conditions that jeopardise their financial sustainability. While mature markets transition to oligopoly structures, the supply chain operation should adapt to a more customer-centric focus. Key issues related to the modelling and impact of the related contractual agreements between firms and customers remain largely unexplored. In the present work, we examine the problem of fair customer allocation in oligopolies under different contractual agreements within a multi-period setting. We consider an ensemble of contract types that vary in terms of pricing mechanisms and duration. The role of fairness is examined following the social welfare and Nash bargaining scheme. In the latter case, the overall problem is formulated as an MINLP. For its efficient solution we employ a piecewise linearisation strategy based on special-ordered sets. The impact of the different fairness schemes on the optimal customer allocation is evaluated via two case studies from the  industrial gases market.
\end{abstract}
\end{frontmatter}


\section*{List of symbols}
\subsection*{Sets}
\begin{tabbing}
$b$   \hspace{5em} \=   Outsourcing tiers\\
$c$     \>   Customers\\
$C_f$       \> Existing customers of firm f\\
$cti$ \>Set of customer's tanks for product i\\
$f$       \>Oligopoly firms\\
$i$         \>Liquid products\\
$j$          \> Liquid and gas products\\
$k$       \>Contracts\\
$n$         \> Grid points\\
$p$      \>Time periods\\
$pe$      \> Set of swap balancing intervals\\
$t$             \> Customer tanks\\
$s$          \> Contract terms
\end{tabbing}   

\subsection*{Parameters}
\begin{tabbing}
$\alpha_f$ \hspace{5em} \= Negotiation power of firm f\\
$\beta_{ictfk}$ \>Base price for selling product i to customer c and tank t by firm f under contract k\\
$\varepsilon$ \>Acceptable power fluctuations from contracted energy consumption\\
$\epsilon_{fkps}$   \> Escalation factor of firm f for contract k in time period p and contract term s\\
$\eta_{f'f}$ \>  Inter-firm swaps premium\\
$\zeta_{b}$ \>  Outsourcing premium\\
$\xi$ \>  Inter-firm swap bound\\
$\pi_f^{sq}$ \> Status quo profit of firm f prior to the fair allocation of the customers (\$)\\
$\pi_{fn}$ \> Profit of firm f at grid point n (\$)\\
$a_{ifp}^{L}$    \> Lower bounding parameter of product i inventory of firm f in time period p\\
$a_{ifp}^{U}$    \> Lower bounding parameter of product i inventory of firm f in time period p\\
$CC_c^{max}$ \>Maximum purchase cost of customer c (\$)\\ 
$CC_{cn}$ \>  Purchase cost of customer c at grid point n (\$)\\
$CEC$ \> Contracted energy consumption (GWh)\\
$D_{ictp}$ \>Demand of product i by customer c in tank t and time period p ($m^3$)\\
$DC_{ictfp}$ \>Delivery cost of product i for customer c in tank t  served by firm f in time period p ($m^3$)\\
$E_{cfk}$ \> Existing customers c of firm f with contract k\\
$EP_{p}$ \> Expected electricity price in time period p (\$/MWh)\\
$FDC_{cf}$ \> Fixed cost of firm f for dropping customer c (\$)\\
$FNC_{cf}$ \> Fixed cost of firm f for acquiring customer c (\$)\\
$INV^0_{if}$ \> Initial inventory for product i in firm f (m\textsuperscript{3})\\
$INVC_{if}$\> Unit inventory cost of product i by firm f (\$/m\textsuperscript{3})\\
$L_k$      \> Duration of contract k\\
$OC_{ictfb}$ \> Piecewise constant outsourcing premium cost of tier b (\$/m\textsuperscript{3})\\
$OT$ \> Average operating hours of the ASU unit (h)\\
$P_{ictfkp}$ \> Price of product i for customer c and tank t served by firm f in time period p\\
\>and contract k (\$/m\textsuperscript{3}) \\
$SWC_{ictff'p}$ \> Swapping premium unit cost\\
$T_{cs}$ \> Formula contract term for customer c\\
$UPC_{if}$  \>Unit production cost of product i by firm f (\$/m\textsuperscript{3})\\
$USC_{ictfp}$ \> Unit service cost of demand of product i for customer c and tank t served by firm f \\
\> in time period p and contract k (\$/m\textsuperscript{3})\\ 
$VDC_{cf}$ \> Variable cost of firm f for dropping customer c (\$)\\
$VNC_{cf}$ \> Variable cost of firm f for acquiring customer c (\$)
\end{tabbing}

\subsection*{Binary Variables}
\begin{tabbing}
$W_{cfkp}$ \hspace{3em} \= 1, if a customer c is served by a company f with contract k in time period p\\
$WD_{cfp}$ \> 1, if customer c is dropped by a company f under contract k in time period p\\
$WN_{cfp}$ \> 1, if customer c is acquired by a company f under contract k in time period p\\
$WS_{cfkp}$ \> 1, if customer c initiates a contract k with firm f in time period p\\
$Y_{ifbp}$ \> 1, if the production of product i is allocated by firm f on the spot market in time period p
\end{tabbing}

\subsection*{ Continuous Variables}
\begin{tabbing}
$\hat{O}_{ictfpb}$ \hspace{3em} \= Disaggregated level of outsourcing for demand of product i for customer c and tank t\\ \>served by firm f in time period p (m\textsuperscript{3})\\
$\delta_{fp}^+$\> Energy consumption deviation above the contracted energy limits (GWh)\\
$\delta_{fp}^-$\> Energy consumption deviation below the contracted energy limits (GWh)\\
$\theta_{fp}$\> Energy consumption within the contracted energy limits (GWh)\\
$\lambda_{fn}$ \> SOS2 variables associated with the piecewise linear approximation\\
\>of the profit of firm f in grid point n \\
$\mu_{cn}$ \> SOS2 variables associated with the piecewise linear approximation\\
\>of the purchase cost of customer c in grid point n \\
$\pi_{fp}$ \> Profit of firm f in time period p (\$)\\
$CAP_{ifp}$ \> Firm's capacity for product i in time period p (m\textsuperscript{3})\\
$CC_{c}$  \> Purchasing cost of customer c in time period p (\$)\\
$EC_{fp}$ \> Electricity cost of firm f for time period p (\$)\\
$IC_{fp}$ \> Inventory cost of firm f for time period p (\$)\\
$INV_{ifp}$ \> Product i inventory of firm f in time period p (m\textsuperscript{3})\\
$NC_{fp}$  \> Customer acquisition cost of firm f for time period p (\$)\\
$O_{ictfp}$ \> Outsourcing amount of product i to customer c and tank t\\
\>by firm f for time period p (m\textsuperscript{3})\\
$PW_{fp}$ \> Power consumption of firm f in time period p (kW)\\
$Q_{ifp}$ \> Firm's in house production of liquid product i in time period p (m\textsuperscript{3}) \\
$RC_{fp}$  \> Customer forfeit cost of firm f for time period p (\$) \\
$S_{ictfp}$ \> Amount of product i sold to customer c and tank t\\
\> from in-house production of firm f for time period p (m\textsuperscript{3})\\
$SC_{fp}$  \> Service cost of firm f for time period p (\$)\\
$SW_{ictf'fp}$ \> Swap amount of product i to customer c and tank t \\ \>from firm f' to firm f in time period p (m\textsuperscript{3})\\
$V^{air}_{fp}$ \> Volumetric flow rate of air in the ASU of firm f in time period p (m\textsuperscript{3}/h) \\
$V_{jfp}$      \> Volumetric flow rate of products j in the ASU of firm f in time period p(m\textsuperscript{3}/h)\\
$V_{ifp}$  \>           Volumetric flow rate of liquid product i in the ASU of firm f in time period p (m\textsuperscript{3}/h)\\
\end{tabbing}

\section{Introduction}
In recent years, the volatile and capital intensive environment has impacted the market structure of process industries. The need for cost minimisation and resilience against supply chain disruptions, has led to the accumulation of production and distribution networks to a limited number of companies, thus creating oligopolies in various sectors. Examples can be found in the steel, oil and industrial gases industries, where a small number of firms accommodates the demand offering similar products. Studying the market dynamics in oligopolies is of high importance to safeguard the normal operation of the supply chains, both in terms of customer demand satisfaction and the aversion of monopoly transformation. Even though the aforementioned goals are traditionally satisfied when the higher level of competition between firms is considered, in this study we investigate an alternative approach based on cooperation. In the proposed framework, the allocation of the customers under different cooperative fairness schemes is evaluated aiming to maintain the market structure while allowing a profit increase for the stakeholders.

Game theory has been extensively studied within the process systems engineering area, a recent review paper on the field can be found in \cite{Marousi2023}. The design of supply chains and operations management  (\cite{leng2005, nagarajan2008,omegaLR}) have been widely facilitated from the use of a game theoretic approach. Games can be divided into categories depending on the player interaction. Those can be cooperative games, were players are collaborating to achieve a win-win outcome, or competitive games were each player acts individually aiming to maximise their individual profit. In the latter case, competitive games can be either Cournot \citep{cournot1838}, in which players make decisions simultaneously, or Stackelberg games \citep{Stackelberg} in which the leader makes a decision first followed by the follower's decision. The choice of game structure impacts the model formulation and consequently different suitable solution approaches. 

\cite{Gjerdrum2001} utilised a game theoretic approach to solve the fair profit allocation in a multi-enterprise supply chain. In the following year, the same authors \citep{gjerdrum2002fair} evaluated two different solution approaches, spatial Branch and Bound and McCormick relaxations, to solve the problem of fair transfer price and inventory optimisation. \cite{Zhao2010} investigated a decentralised supply chain model between manufacturers and suppliers to determine the optimal wholesale contract selection. The authors found that using a Nash bargaining model allowed both parties to maximize their profit under the chosen wholesale price mechanism. \cite{yue2014fair} proposed a logarithmic transformation and Branch and Refine (BR) algorithm to solve the optimal operation and profit allocation of bio-ethanol supply chain under a Nash bargaining scheme. On the other hand, \cite{Zheng2019} utilised post-optimisation fairness criteria, i.e. Shapley and nucleolus values, for the profit allocation in a three-echelon closed loop supply chain and equal satisfaction mechanisms to evaluate fairness and find the optimal profit allocation in the examined supply chain. The impact of different fairness schemes for payoff allocation have been examined in liquid gas supply chains \citep{Charitopoulos2020},eco-industrial park \citep{cruz2021} and multi-carrier energy system \citep{omegafair2} design problems among others. \cite{omegafair1} observed that the allocation based on Shapley value, hence the dissatisfaction of the coalitions, could be powerless in logistic problems with incomplete information, and thus proposed an alternative method  that reflects the dissatisfaction of the players instead. The selection of a fairness approach is not always straightforward, given that there may be conflicting objectives such as profit profit maximisation while satisfying sustainability commitments. \cite{Koleva2018} have evaluated such a conflicting problem in order to find the economically and environmentally optimal design of a water supply chain under uncertainty. 

From a non-cooperative perspective, \cite{Levis2007} introduced a Bertrand-type model to determine Nash equilibrium prices for maximising profits in a duopoly. \cite{Zamarripa2012} assessed the impact of different degree of cooperation, both cooperative and non-zero-sum games, in a supply chain optimisation problem. To represent the hierarchical structure of supply chain problems, multi-level programming and decentralized game-theoretic approaches are often proposed.\cite{Yue2014} modeled a biofuel supply chain as a Stackelberg game using bi-level programming. The authors employed an improved Branch and Refine algorithm to retrieve the global optimal solution by solving a series of MILP sub-problems. The same authors \citep{Yue2017} further investigated reformulation and decomposition algorithms for bi-level MILP problems stemming from Stackelberg game formulations. Multi-objective formulations arise in non-cooperative games as well. For example, \cite{Gao2017} considered a non-cooperative shale gas supply chain as a multi-objective bi-level MILP problem, aiming to maximise the net present value and minimise the greenhouse gas emissions. 

Even though in this study the contract formulations and contract terms are parameters of the model and do not constitute decision variables, a brief review on the contract selection is provided for completeness. \cite{Park2006} have addressed the problem of contract modelling in a multi-period framework, including contract selection in supply chain models. Disjunctive programming was employed to address both long term and short term operations. The contracts proposed for supply/demand are: a) fixed price, b) discount after certain amount, c) bulk discount, and d) fixed duration. A classification of different supply contract in a multi-period programming problem for optimal contract selection was examined by \cite{Bansal2007}. The study incorporated the dominant real-life contract features, such as purchase commitments and flexibility, commitment duration and bulk prices/ discounts.
\cite{Qin2007} evaluate a non-cooperative Stackelberg game, where supplier acts as the leader and decides on a pricing policy, the buyer reacts a follower and determines the annual sales volume. For this application volume discounts are considered. \cite{Calfa2015} have incorporated the optimal contract selection in the scheduling problem of a chemical process network, by choosing among the contracts proposed by \cite{Park2006}.\cite{Qi2016} studied the scheduling problem of energy intensive process industries by proposing a block contract model so as to take into account widely used power contracts. In order to propose an economically viable integration of renewable energy generators in the existing energy markets, \cite{omegacontracts} used a game-theoretic using physical and financial contracts as a mean of evaluating the different market structures. Discount contracts for supplier/ manufacturer agreements in the process industry can be  found in the papers of \cite{martin2018, Kirschstein2019}. Recently, \cite{Elekidis2022} examined the optimal contract selection for contract manufacturing organisations in the in the secondary pharmaceutical industry, considering demand uncertainty and risk measures such as Value-at-Risk and Conditional Value-at-Risk. 

The present study focuses on the fair customer allocation of liquid gas market supply chains. A game-theoretic framework is proposed to examine customer re-allocation for profit maximisation of oligopoly firms. In the status quo market there is a set of free customers that are not supplied by any firm and allow for a market share growth. The introduced models transform the static model of \cite{Charitopoulos2020} in a multi-period framework. Under this scope, customers can alternate between different firms and contractual agreements in order to maximise the firms' profit. Note that the price of the products and the contract formulations are parameters of the model. The aim of this paper is to (i) evaluate the impact of contract choice in customer allocation and process operation in oligopolies and (ii) apply different fairness schemes so as to distribute the profits of the market expansion among the oligopoly firms. The paper has the following structure:the problem statement is introduced in Sect.\ref{sec:Prob} and the multi-period model in Sect.\ref{sec:Model}.The two case studies are examined in Sect.\ref{CaseGT} along with the corresponding results. Finally, the key findings and future directions are summarised in Sect.\ref{sec:conc}.

\section{Problem statement}\label{sec:Prob}
Recently, a static game-theoretic approach for the fair customer allocation within oligopolies was proposed by \cite{Charitopoulos2020}. In this work a multi-period framework is proposed, where we investigate the role of fairness for the tactical allocation of customers within existing oligopolies under different contract options. In the present work we take the viewpoint of the oligopoly firms  which seek to maximise their profits in a fair manner. To this end, we assume that the firms that constitute the oligopoly are rational and that each firm has estimates of the other firms' information. Customers that are served by the oligopoly have to be assigned to one of the firms and at the beginning of the planning horizon are split between existing customers that hold a contract with one firm and new customers that provide an opportunity for market share growth. The overall problem is formulated as a multi-period mixed integer nonlinear program (MINLP), which is reformulated to an MILP class,  whilst the role of fairness is examined via the social welfare and Nash Bargaining approaches. Overall the problem statement has as follow:

\hspace{-0.6cm}\textbf{Given:}
\begin{itemize}
\item Tactical horizon for supply chain planning
\item A portfolio of existing and new customers
    \item Customer demand per liquid product
    \item Firms' plants production and inventory capacities
    \item A set of different contract types that differ in duration \& customer pricing
    \item Firms' target levels on customers' contract mix   
    \item Product prices per customer and firm
    \item Delivery cost to the customer's tank
    \item Customer acquisition/forfeit variable and fixed costs
    \item Third-party production costs and tiers
    \item Customer demand "swap" cost and limits between firms
    \item Electricity  cost projection and consumption agreement with energy system operator (ESO)
\end{itemize}
\textbf{Compute:}
\begin{itemize}
    \item Optimal customer assignment to the firms
    \item Contract selection and schedule between firms and customers
    \item Optimal production \& inventory profiles for the firms
    \item Optimal product demand swap levels for the firms
    \item Optimal product demand outsourcing levels for the firms

\end{itemize}
\textbf{So as to:}
\begin{itemize}
    \item \textit{Fairly} maximise the firms' profit 
\end{itemize}

\begin{figure}[]
\centering
\includegraphics[width=12cm]{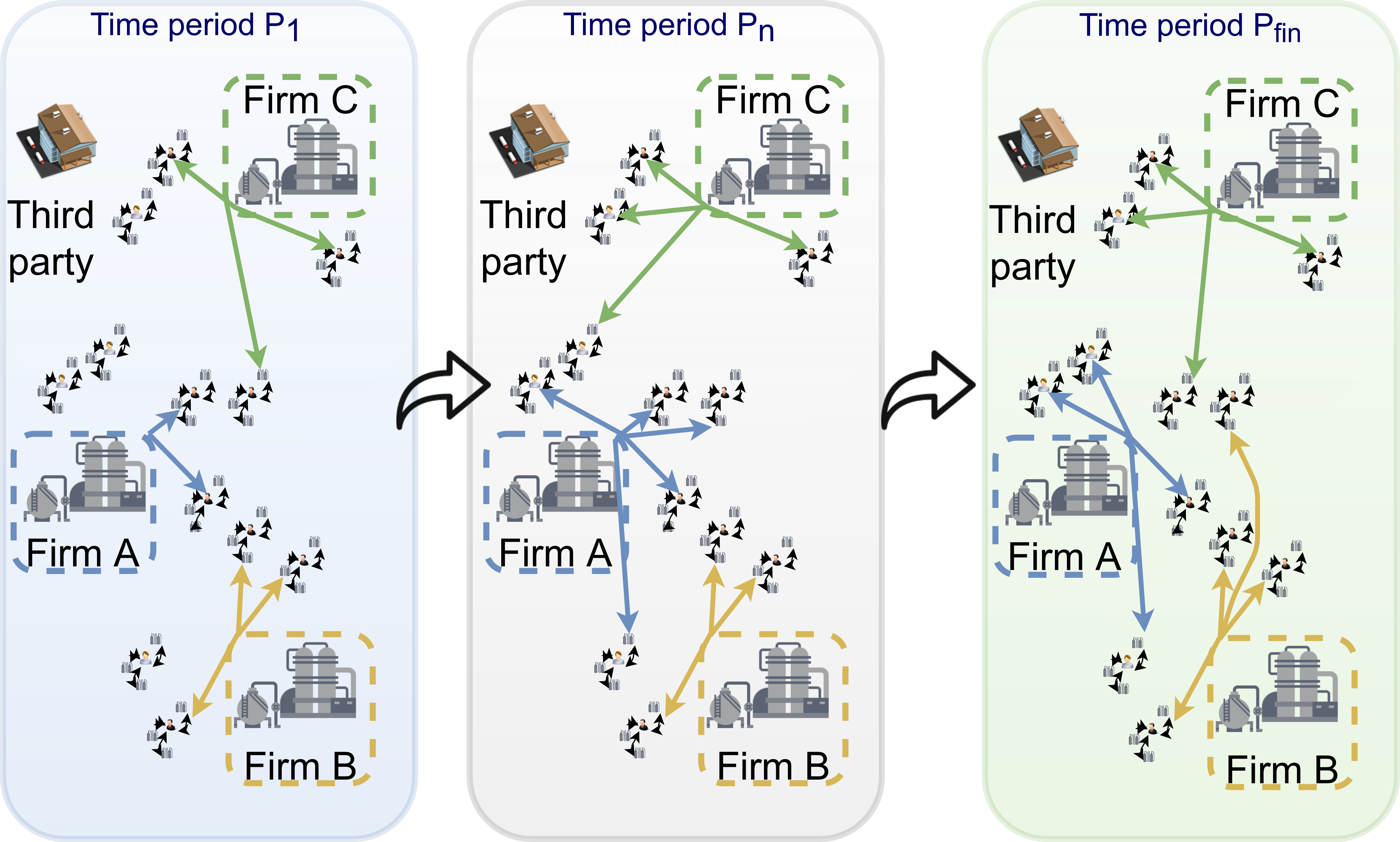}
\caption{Conceptual representation of an oligopoly comprised by three companies that serve a number of customer tanks.}\label{fig:oligopoly}
\end{figure} 

A conceptual representation of the problem under study is shown in Fig. \ref{fig:oligopoly}. At the beginning of the planning horizon each firm has a set of assigned customers and there exist new customers that have not been assigned to any company. In the subsequent time periods, depending on the contractual agreement between the firms and the customers, a customer may be assigned to another firm with a new contract or remain with the same firm with the same or different contract. The mathematical formulation for the multiperiod model is presented next in Section \ref{sec:Model}.

\section{Multi-period model}\label{sec:Model}
In this section, first the model formulation for the customer allocation in oligopolies is presented and next, the game-theoretic framework for its fair optimisation is introduced. The nomenclature of the mathematical developments is provided in the beginning of the article. The key assumptions in the present work are summarised as follows: (i) firms will participate in the game only if they can achieve greater profit than their current one, (ii) deterministic production and service cost, (iii) customer demands are given as their deterministic average values, (iv) firms should serve all customers and (v) decentralised decision making.

\subsection{Customer assignment and contracts scheduling}
At any time period ($p$) customers ($c$) must be served by one firm ($f$) and under one type of contract ($k$). This condition is modelled by introducing the binary variable $W_{cfkp}$ in  Eq.\eqref{eq:binassign}. The starting of a new contract is recorded by the binary variable $WS_{cfkp}$ as modeled in Eq.\eqref{eq:startup}. It is noteworthy that all customers are assigned in the first time period ($p=1$) and no new customers are introduced in later time periods.
\begin{align}
    & \underset{f}{\sum}\underset{k}{\sum}W_{cfkp}=1 \qquad \forall c,p \label{eq:binassign}
\end{align}
\begin{align}
    & W_{cfkp}-W_{cfk,p-1} \leq WS_{cfkp} \qquad \forall f,c,k,p \label{eq:startup}
\end{align}
Furthermore, at any time period customers may sign at most one contract with one firm. To model this instance, the binary variable $WS_{cfkp}$ is introduced which denotes whether at time period $p$ customer $c$ signed contract $k$ with firm $f$ along with Eq.\eqref{eq:StartContract}.

\begin{align}
    & \underset{f}{\sum}\underset{k}{\sum}WS_{cfkp}\leq 1 \qquad \forall c,p \label{eq:StartContract}
\end{align}

\begin{figure}
    \centering
    \includegraphics[width=0.8\linewidth]{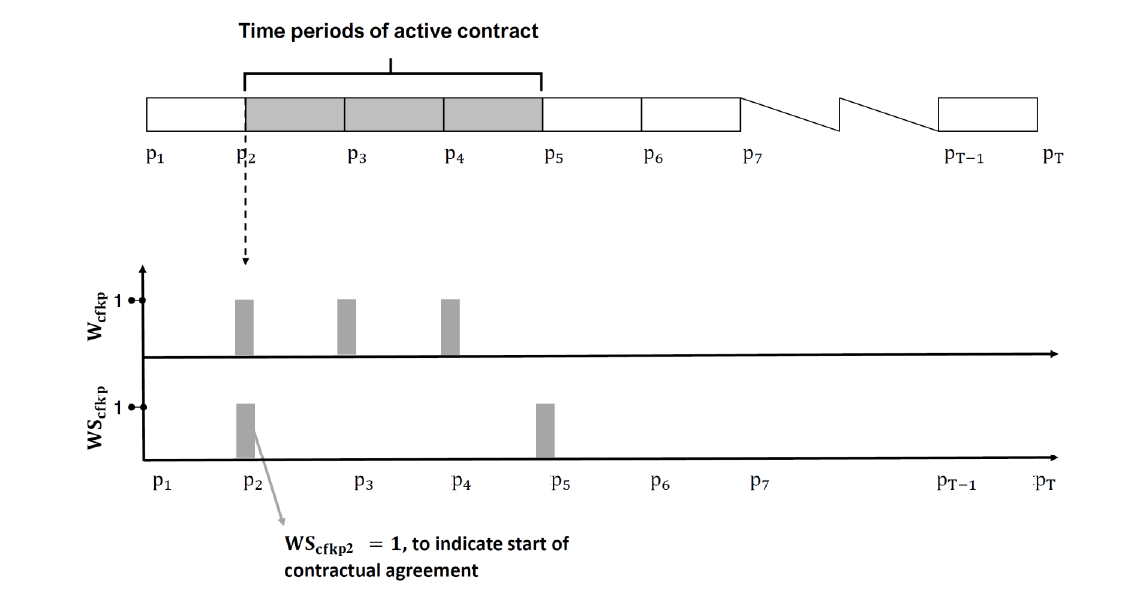}
    \caption{Relation between the binary variables $WS_{cfkp}$ which denote the initiation of a contract and the binary variables $W_{cfkp}$ which denote that a contract is active over the periods spanning its agreed duration $L_{k}$.}
    \label{fig:contracts}
\end{figure}

\noindent The modeled contracts are closed, i.e., once a contract is signed it must be enforced during its duration $L_k$. To eliminate the possibility of a customer signing a contract prior to the end of their current one Eq.\eqref{eq:mintime} is introduced. A conceptual representation of Eq.\eqref{eq:mintime} is given in Fig.\ref{fig:contracts}, in time period $p=4$, a  3 period contract is active, thus $W_{cfk,p=4}=1$, this means that from $p-L_k +1= 4-3+1= 2^{nd}$ period to the $4^{th}$ period only one new contract $WS_{cfk}$ can be signed.
\begin{equation}
\sum_{p-L_{k}+1}^{p}WS_{cfkp'}=  W_{cfkp} \qquad \forall c,f,k, p\geq L_{k} \label{eq:mintime}
\end{equation}
The acquisition and forfeiting of a customer by a firm is modeled by Eqs.\eqref{eq:newbin} and \eqref{eq:forfbin}, where the binary variables $WC_{cfp}$ and $WD_{cfp}$ mark  respectively the acquisition and forfeiting of a customer by a firm in a specific time period.
\begin{align}\label{eq:newbin}
 \sum_k WS_{cfkp}- WS_{cfk,p-L_k}\leq WC_{cfp}    \qquad \forall c,f, p\geq 1
\end{align}

\begin{align}\label{eq:forfbin}
 \sum_k WS_{cfk,p-L_k}- WS_{cfkp} \leq WD_{cfp}    \qquad \forall c,f, p\geq 1
\end{align}

\subsection{Contract formulation}\label{sec:contacts}
For the examined case studies, different "Formula" contracts have been evaluated.  The selling price of a product unit is determined by taking into account the base parameter ($\beta_{ictfk}$), a number of key parameters ($T_{cs}$) which are escalated by the escalation factor ($\epsilon_{fkps}$). Representative values for the contract parameters and escalation factors can be found in the appendix in Fig. \ref{fig:Terms} and \ref{fig:escal} respectively. The pricing of the contracts is a deterministic parameter in the examined models. For the starting period of the model, Eq.\eqref{eq:FormIn} dictates the price that depends on the base and the key parameters of the contract, as time progresses ($p>1$) the escalation factor is introduced so as to introduce a time variance to the product price (Eq.\eqref{eq:Formula}).
\begin{equation}
    P_{ictfkp}=\beta_{ictfk}\sum_s T_{cs}, \qquad p=1, \forall i,c,t,f,k \label{eq:FormIn}
\end{equation}
\begin{equation}
    P_{ictfkp}=\beta_{ictfk}\big(\sum_s T_{cs}(1+\epsilon_{fkps})\big), \qquad p>1,~\forall i,c,t,f,k \label{eq:Formula}
\end{equation}
\subsection{Customer demand satisfaction}\label{sec:demand}
\begin{figure}
\centering
\includegraphics[width=1\linewidth]{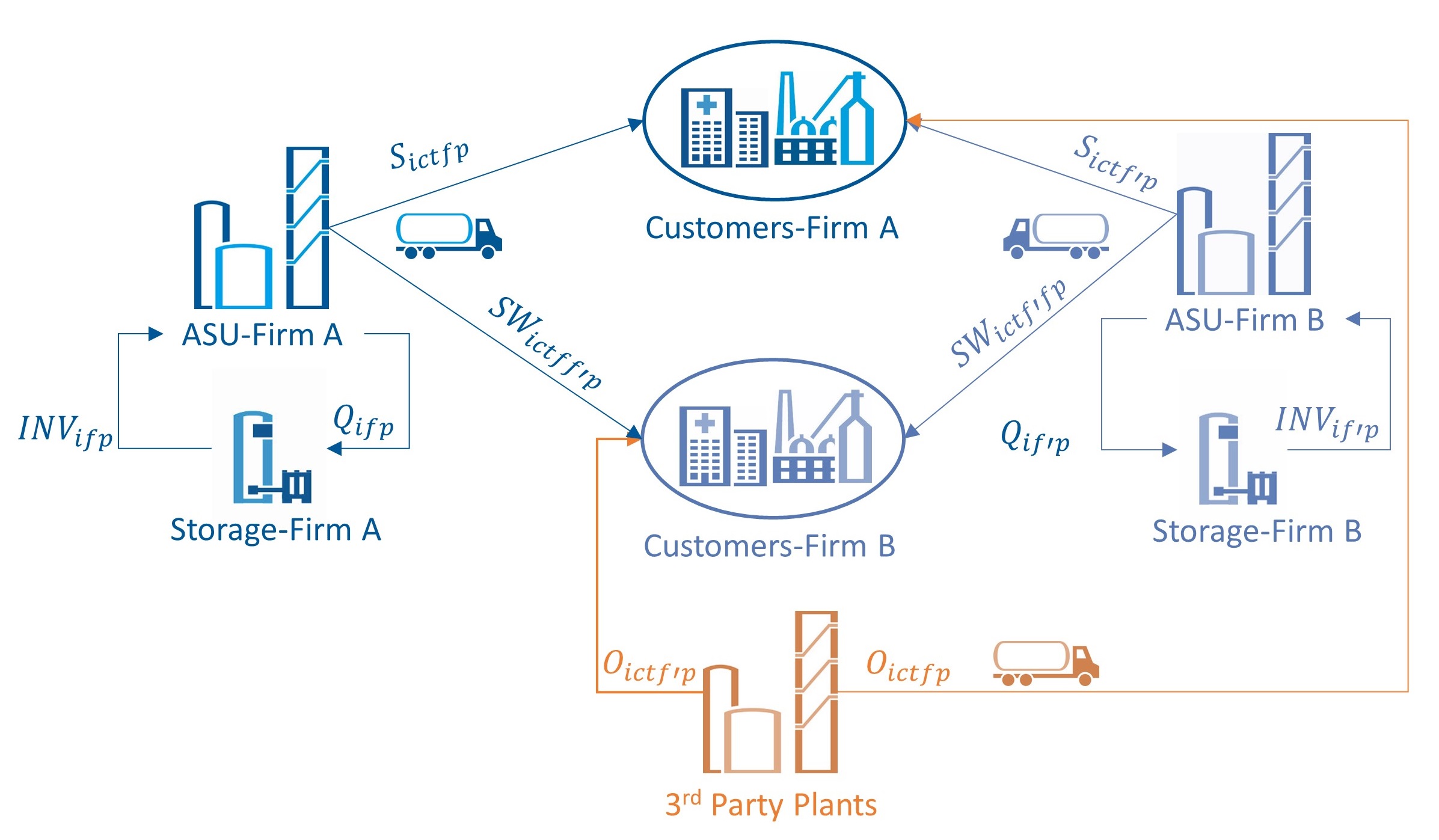}
\caption{\label{fig:Demand-satisfaction-mechanisms} Product flows and customer demand service in an duopoly with a spot market and inventory sites.}
\end{figure}

Each customer has a set of tanks dedicated to particular products. Once a customer has been contracted to a firm, then the firm should satisfy their demand for specific products. This can be either done via product amount covered by the contracted firm's in-house production ($S_{ictfp}$), product swaps with other firms ($SW_{ictf'fp}$) or via purchasing product amounts from the spot market ($O_{ictfp}$). To model this Eq.\eqref{eq:demand_satisfaction} is employed, and the aforementioned concepts are illustrated in Fig.\ref{fig:Demand-satisfaction-mechanisms}. Each firm is represented by a single Air Separation Unit (ASU), the capacity of which provides an upper bound to the quantity of product produced to supply the assigned customers and the exchanges with the other firms of the oligopoly, Eq.\eqref{eq:capacity}.

\begin{align}
& S_{ictfp}+\underset{f'\neq f}{\sum}SW_{ictf'fp}+O_{ictfp}=D_{ictp}\sum_{k}W_{cfkp}\qquad\forall i,c,t\in cti, f, p\label{eq:demand_satisfaction} 
\\
& Q_{ictfp}+\underset{f'\neq f}{\sum}SW_{ictff'p}\leq Cap_{ifp}\qquad\forall i,c,t\in cti, f, p\label{eq:capacity}
\end{align}

Note that the in-house capacity for each of the products ($Cap_{ifp}$) is equal to the volumetric flow rate of each liquid product for the examined time period adjusted accordingly to the time resolution of the model. Finally, to capture the inventory dynamics of the supply chain Eq. \eqref{eq:inventorybalance} is employed along with Eq.\eqref{eq:inventorycapacity} to model inventory storage capacity bounds. 
\begin{multline}
INV_{ifp} = INV_{if,p-1_{|p>1}} + INV^0_{if,p=1} + Q_{ifp} \\-\sum_{(c,t) \in cti} \sum_{f'\neq f}SW_{ictff'p} - \sum_{(c,t)\in cti}S_{ictfp} \qquad \forall i,f,p \label{eq:inventorybalance} 
\end{multline}

\begin{align}
& a^{L}_{ifp}V_{ifp} \leq INV_{ifp} \leq  a^{U}_{ifp}V_{ifp}  \qquad \forall i,f,p\label{eq:inventorycapacity}
\end{align}

\subsection{Plant production surrogate model}
Air separation plants have high energy requirements in order to convert atmospheric air to different gaseous and liquid products, i.e. gaseous nitrogen ($\mathrm{GNI}$), gaseous oxygen ($\mathrm{GOX}$), liquid nitrogen ($\mathrm{LNI}$), oxygen ($\mathrm{LOX}$) and argon ($\mathrm{LAR}$). Let $J$ be the set of the products of an ASU unit so as $J=:\{LOX,LNI,LAR,GOX,GNI\}$.  While this work primarily focuses on liquid products, these are secondary products in the gas market, since air separation units (ASU) are designed to produce gaseous products. It is important to recognize that the underlying characteristics of the ASU directly impact the firm's production capacity in meeting customer demand. To address these considerations regarding production capacity and electricity consumption, a surrogate model is employed. This surrogate model effectively captures the correlation involved in manufacturing different liquid products and the corresponding energy consumption. The general form of the model is given by Eq.\eqref{shortcut_1}-\eqref{shortcut_2}. An outline of the equations of the surrogate model can be found in \cite{Charitopoulos2020}.

\begin{align}
    & PW_{fp}= f(V_{fp}^{air},V_{jfp}) \qquad j \in J, \quad\forall f,p \label{shortcut_1}\\
    & g(V_{fp}^{air}, V_{jfp})\leq 0 \qquad \qquad j \in J,\quad\forall f,p \label{shortcut_2}
\end{align}

\noindent where ${PW_{fp}}$ is the electricity power consumed and ${V_{fp}^{air}, V_{jfp}}$ stand for the volumetric flows of the air and the different gaseous and liquid products.

\subsection{Spot market product acquisition}
To account for instances where production and inventory capacities do not suffice to meet customer demand, product acquisition from the spot market at typically higher cost. Spot market product acquisition quantities are represented using the positive variable $O_{ictfp}$. The cost of spot market acquisition is assumed to follow a piecewise constant paradigm, to model this, $b$, not necessarily continuous, tiers are employed. Each tier inflicts specific spot market premium costs ($OC_{icftb}$) as shown by Eqs.\eqref{eq:outsource1}-\eqref{eq:outsource4}. The positive variable  $\hat{O}_{ictfbp}$ is employed to model the disaggregated counterpart, per tier, of  $O_{ictfp}$ while the binary variable $Y_{ifbp}$ is active for only one tier each time period when product is purchased from the spot market. 

\begin{align}
 & O_{ictfp}=\underset{b}{\sum}\hat{O}_{ictfbp}\qquad\forall i,c,t \in cti ,f,p\label{eq:outsource1}\\
 & \gamma_{b}^{L}Y_{ifbp}\leq\underset{c,t\in cti}{\sum}\hat{O}_{ictfbp}\leq\gamma_{b}^{U}Y_{ifbp}\qquad\forall i,f,b,p\label{eq:outsource2}\\
 & \underset{b}{\sum}\hat{O}_{ictfbp}\leq D_{ictp}\sum_{k}W_{cfkp}\qquad\forall i,c,t \in cti,f,p\label{eq:outsource3}\\
 & \underset{b}{\sum}Y_{ifbp}\leq1\qquad\forall i,f,p\label{eq:outsource4}
\end{align}
Eq.\eqref{eq:outsource1} is employed to disaggregate the overall outsource amount into different tiers and the lower and upper bounds of each tier are given by Eq.\eqref{eq:outsource2}. Eq.\eqref{eq:outsource3} is a logic constraint which implies that customer demand can be acquired by the spot market only if the customer is currently served by the firm and finally Eq.\eqref{eq:outsource4} defines that at most one tier can be selected for spot market acquisition  of product demand for each product type ($i$) and firm ($f$).  

\subsection{Inter-firm swap agreements}
As mentioned in Section \ref{sec:demand}, a firm can hold swap agreements with another firm so as to achieve a lower service cost in case a customer is closer to the not contracted firm. Compared to spot market product acquisition, inter-firm swaps rely on bilateral contracts, which specify terms and conditions on product amounts that can be swapped. To model this, Eqs. \eqref{eq:swap_bound}-\eqref{swap_logic}
are stipulated, which express that the total amount of product demand that can swapped  ($SF_{ictff'p}$) within a time period should be less than a pre-specified upper bound ($U_{if}$) and firms may only employ swaps for customers within their portfolio. 
\begin{align}
&\underset{c,t\in cti,f'\neq f}{\sum}SW_{ictff'p}\leq \xi Cap_{ifp} \qquad\forall i,f,p \label{eq:swap_bound} \\
&SW_{ictff'p}\leq D_{ictp}\sum_{k}W_{cfkp} \qquad\forall i,c,t \in cti, f,f'\neq f,p \label{swap_logic}
\end{align}

In order to avoid the transformation of the oligopoly in a monopoly, where basically only one firm would be the producer and the other firms would be transformed in retailers, the amount of each product swapped between firms is balanced for a decided time interval ($pe$), Eq.\eqref{eq:swap_condition}. At the same time, such a constraint prevents any convolutions stemming from the pricing of swapped product between firms.

\begin{align}
& \underset{p \in pe}{\sum}~\underset{c,t \in cti}{\sum}SW_{ictff'p}=\underset{p \in pe}{\sum}~\underset{c,t \in cti}{\sum}SW_{ictf'fp}\qquad\forall i, f, f'\neq f \label{eq:swap_condition}
\end{align}

\subsection{Customer service cost}

The customer service cost  (${SC_{fp}}$) accounts for the different demand satisfaction mechanisms under consideration. Firstly the unit service cost of serving a specific customer tank (${USC_{icfp}}$) is calculated as ${USC_{ictfp}=\frac{DC_{ictfp}}{D_{ictp}}}$, where ${DC_{ictfp}}$ is the delivery cost of firm f for serving customer's c tank t. The first term as shown by Eq.\eqref{eq: serving cost} reflects the monthly average delivery cost for the case of in-house production. The second term reflects the cost incurred by the swaps, where a swap premium is involved
(${SWC_{iff'p}}$) and is assumed to involve $\eta_{f'f}$times higher costs, depending on the inter-firms contracts established, for the firm compared to the in-house production and delivery as shown by Eq.\eqref{eq:swapping premium}. The last term of Eq.\eqref{eq: serving cost} reflects the cost of spot market product acquisition, where similarly to the swap's cost, a spot market premium is considered depending on the related tier. Each tier is associated with a different premium ($\zeta_{b}$).
\begin{multline}
     SC_{fp}=\underset{i,c,t \in cit}{\sum}USC_{ictfp}S_{ictfp} +\underset{i,c,t \in cit,f'\neq f}{\sum}SWC_{icf'fp}SW_{ictf'fp} \\+ \underset{i,c,t \in cit,b}{\sum}OC_{icfbp}\hat{O}_{ictfpb}\qquad\forall f\label{eq: serving cost}    
\end{multline}

\begin{align}
& SWC_{icf'fp}=\eta_{f'f} USC_{icfp}\qquad\forall i,c,f'\neq f, p \label{eq:swapping premium}
\\
& OC_{icfbp}=\zeta_{b}(USC_{icfp}+UPC_{ip})\qquad\forall i,c,f,b, p \label{eq:outsourcing premium}
\end{align}

\subsection{Customer acquisition cost}
At the beginning of the planning horizon there exist a number of customers that have not been contracted to any firm which provide opportunity for market share growth. Apart from these customers, a firm may acquire a customer that was previously assigned to different firm. In any of these cases, customer acquisition inflicts a set of costs related to administrative duties as well as installation and maintenance of product tanks on the customer's site. To model these costs, a fixed ($FNC_{cf}$) and a variable cost component ($VNC_{ctf}$) are considered, with the variable cost being proportional to customer's demand. To accurately track changes in customer assignment to firms, the binary variable $WN_{cfp}$ is introduced which denotes new customer acquisition as shown by Eq. \eqref{eq:newcustomerbinary}. Analogously, the binary variable corresponding to the forfeit of a customer by a firm $WD_{cfp}$ is modeled in Eq.\eqref{eq:dropcustomerbinary}
\begin{align}
        & WN_{cfp} \geq \sum_{k}(WS_{cfkp}- WS_{cfk,p-L_{k}}) \qquad \forall c,f,k,p>L_k \label{eq:newcustomerbinary}
\end{align}
\begin{align}
    & WD_{cfp} \geq \sum_{k}(WS_{cfk,p-L_{k}} - WS_{cfkp})  \qquad \forall c,f,k,pp>L_k\label{eq:dropcustomerbinary}
\end{align}

The new customer acquisition cost can be computed by Eq.\eqref{eq:acq_cost_p1}-\eqref{eq:acq_cost}.
\begin{multline}
NC_{fp}=\underset{c\notin C_{f}}\sum(FNC_{cf} + \sum_{(i,t)\in cti}VNC_{ctf}D_{ict})WS_{cfkp} \\\qquad \forall f,p=1\label{eq:acq_cost_p1}    
\end{multline}
\begin{multline}
NC_{fp}=\sum_{c}(FNC_{cf} + \sum_{(i,t)\in cti}VNC_{ctf}D_{ict})WN_{cfp} \\ \qquad \forall f,p>1\label{eq:acq_cost}    
\end{multline}

\noindent Notice that for the first time period of the planning horizon, Eq.\eqref{eq:acq_cost_p1}, considers only new customers using the complement of the set $C_{f}$ which is the set of firms' existing customers.

\subsection{Customer forfeit cost}
To model instances where a customer is forfeited by a firm we consider a fixed ($FDC_{cf}$) and a variable cost component ($VDC_{ctf}$) which reflect administrative as well as tank decommissioning costs among others. To this end, the variable cost is proportional to the customer demand. We assume that change in contract type does not inflict such cost and thus the overall cost of forfeiting a customers is given by Eqs. \eqref{eq:dropping_cost_p1}-\eqref{eq:dropping_cost}. 

\begin{multline}
RC_{fp}=\underset{c\in C_{f}}\sum(FDC_{cf} + \sum_{(i,t)\in cti}VDC_{ctf}D_{ict})\sum_{k}(E_{cfk}-WS_{cfkp}) \qquad \forall f,p=1\label{eq:dropping_cost_p1}    
\end{multline}

\begin{multline}
RC_{fp}=\underset{c}\sum(FDC_{cf} + \sum_{(i,t)\in cti}VDC_{ctf}D_{ict})WD_{cfp} \\ \qquad \forall f,p>1\label{eq:dropping_cost}    
\end{multline}

Notice that Eq.\eqref{eq:dropping_cost_p1} is considered for the first time period and the binary parameter $E_{cfk}$ denotes that an existing customer ($c\in C_{f}$) is assigned to firm f under contract k. Here we assume that at the beginning of the planning horizon the company must decide on renewing the contract or dropping the customer at the first time period, however this constraint can be easily modified to account for carry-over remaining duration from the previous planning horizon. 

\subsection{Power consumption cost}
The dominant factor in an ASU plant operating costs is electricity consumption. Due to their nature of being highly energy intensive processes, ASU plants are typically part of industrial demand side response schemes. Within those, the firm signs an agreement with the energy system operator (ESO) so as to regulate its consumption within pre-specified limits at the benefit of a favourable energy price. Should however the firm consume power below or above those thresholds significant penalties apply. To monitor this requirement, the variables $\delta^{+}_{fp}$ and $\delta^{-}_{fp}$ are introduced to denote upward and downward deviations from the pre-specified energy consumption range, while the variable $\theta_{fp}$ denotes energy consumption within those limits( Eqs.\eqref{eq:elecbound1}-\eqref{eq: elecbound3}).

\begin{align}
    & \delta^{+}_{fp}\geq PW_{fp}-(1+\varepsilon)CEC \qquad \forall f,p \label{eq:elecbound1}\\
    & \delta^{-}_{fp}\geq (1-\varepsilon)CEC-PW_{fp} \qquad \forall f,p \label{eq:elecbound2}\\
    &  (1-\varepsilon)CEC\leq \theta_{fp} \leq (1+\varepsilon)CEC \qquad \forall f,p \label{eq: elecbound3}
\end{align}
\noindent The corresponding power consumption balance  and the total energy consumption ($PW_{fp}$)is formulated in Eq.\eqref{eq:energy}.

\begin{equation}\label{eq:energy}
    PW_{fp} = \theta_{fp} + \delta^{+}_{fp} + \delta^{-}_{fp} \qquad \forall f,p
\end{equation}
\noindent To calculate the energy consumption costs of the plants Eq.\eqref{eq:electricity cost} is employed. 

\begin{gather}
EC_{fp}=EP_{p}OT\theta_{fp}+ 1.2(EP_{p}OT(\delta^{-}_{fp}+\delta^{+}_{fp}))\qquad\forall f,p\label{eq:electricity cost}
\end{gather}

\noindent where $EP_{p}$ is the average electricity unit price of planning period p , $OT$ is the average operating time of the ASU plant for the specified time resolution. 
\subsection{Inventory cost}

Inventory utilisation inflicts additional costs ($IC_{fp}$) as shown by Eq. \eqref{eq:invetorycost}. 

\begin{align}\label{eq:invetorycost}
    & IC_{fp} = \sum_{i}INVC_{if}INV_{ifp}
\end{align}
\noindent where the parameter $INVC_{if}$ denotes the unit inventory cost for each firm and for different products. 

\subsection{Profit calculation}

The profit of each firm is calculated as the difference between the revenue and the total costs incurred by the customers' activity. The revenue is calculated as the selling price of product i multiplied by the resulting product demand from customers served by each firm. The profit $\pi_{f}$ for each firm is given by Eq.\eqref{eq:Profit}.

\begin{equation}
\pi_{f}=\underset{p}{\sum}\underset{i,c,t\in cit,k}{\sum}P_{ictfpk}\cdot D_{ictp}W_{cfkp} -SC_{fp}-RC_{fp}-NC_{fp}-EC_{fp}-IC_{fp} \qquad \forall f \label{eq:Profit}
\end{equation}

\subsection{Game-theoretic objective formulation}\label{sec:nash}
In the problem under study, the status quo market involves free customers that are not being served by neither of the oligopoly firms, allowing for a market expansion. Hence, the status quo profit of the firms correspond to the cumulative profit over the examined time horizon but having a smaller pool of available customers. In order to allocate the payoffs to the players of the games two fairness schemes are examined, which result in different objective function formulations, i.e. the social welfare scheme and the Nash bargaining. A more detailed analysis of fairness schemes used in the PSE literature can be found in the work of \cite{Marousi2023}.

A firm will agree on entering the game if they can achieve profit greater than the one in the status quo ($\pi^{sq}$) as shown by Eq. \eqref{eq:rational}.

\begin{equation}
\pi_{f} \geq \pi^{sq}_{f} \qquad \forall f \label{eq:rational}
\end{equation}
Consequently, the Nash equilibrium of the game is computed as the maximum value of the Nash product which is given by Eq.\eqref{eq:nash_product}.

\begin{equation}
\Phi = \prod_{f}(\pi_{f}-\pi^{sq}_{f})^{\alpha_{f}} \label{eq:nash_product}
\end{equation}

The parameter $\alpha_{f}$ represents the negotiation power of each firm f and the summation of the negotiation power over all the firms equals to 1. The Nash bargaining objective function is non-convex and nonlinear, i.e. Eq.\eqref{eq:nash_product} which can result in computationally demanding problems. In this study we employ a separable programming approach to approximate the Nash product. Initially, the logarithmic transformation is applied on Eq.\eqref{eq:nash_product} resulting in Eq.\eqref{eq:lognash}.

\begin{equation}
\Psi= ln\Phi = \underset{f}{\sum}\alpha_{f}ln(\pi_{f}-\pi^{sq}_{f}) \label{eq:lognash}
\end{equation}

Eq.\eqref{eq:lognash} is still nonlinear but now the objective function is strictly concave. Hence, it can be approximated via the piecewise linear function $\tilde{\Psi}$ over a number of $n$ grid points as shown by Eq.\eqref{eq:sosnash}.
\begin{equation}
\tilde{\Psi}= \underset{f}{\sum}\underset{n}{\sum}\alpha_{f}ln(\tilde{\pi}_{fn}-\pi^{sq}_{f})\lambda_{fn} \label{eq:sosnash}
\end{equation}
where the parameter $\tilde{\pi}_{fn}$ corresponds to the profit of firm f at grid point n and  $\lambda_{fn}$  is a SOS2  variable which implies that only two adjacent grid points take nonzero values and satisfy Eq.\eqref{eq:SOS2}. The special ordered set (SOS) approximation was introduced by \cite{Beale1976} and later employed by a series of studies in PSE \citep{Gjerdrum2001, Gutierez2015,liu2017fair,Charitopoulos2020}. The notion of SOS approximation is that a concave function is discretised over the number of grid points and the original problem solution is approximated as the convex combination two adjacent grid points over a line segment. Increasing the number of grid points, results in finer discretisation and thus error reduction with the number of grid points increasing until sufficient accuracy is achieved. The profit of each company can be calculated by Eq.\eqref{eq:profapprox}. 

\begin{table}[ht]
\centering
\caption{Proposed mathematical models of different games of firms}
\label{tab:models_gf}
 \resizebox{\columnwidth}{!}{
 \begin{tabular}{c c c c c c} 
 \hline
\textbf{Fairness scheme} & \textbf{Constraints} & \textbf{Objective} & \textbf{Class}& \textbf{Scope} &\textbf{Abbreviation}\\ \hline
Nash bargaining & Eq.\eqref{eq:binassign}-\eqref{eq:rational}& Eq.\eqref{eq:lognash} &MINLP & Decentralised & FNS\\
Nash bargaining & Eq.\eqref{eq:binassign}-\eqref{eq:rational},\eqref{eq:SOS2},\eqref{eq:profapprox}& Eq.\eqref{eq:sosnash}  & MILP&Decentralised & FLNS\\
Social welfare& eq.\eqref{eq:binassign}- \eqref{eq:Profit}& Eq.\eqref{eq:naiveapp} &  MILP &Centralised& FSW\\
 \hline
\end{tabular}} 
\end{table}

\begin{align}
 &\underset{n}{\sum}\lambda_{fn}=1\qquad \forall f \label{eq:SOS2}\\
&\pi_{f} = \underset{n}{\sum}\tilde{\pi}_{fn}\lambda_{fn} \qquad \forall f \label{eq:profapprox}
\end{align}

From an alternative perspective, the objective function based on the social welfare scheme is formulated in Eq.\eqref{eq:naiveapp}. Table \ref{tab:models_gf} summarises the fairness schemes, equations and scope of the proposed games of firms. Solving the nonlinear Nash bargaining model (FNS) was found intractable when utilising the state of the art solvers, hence this study focuses on the results of the linearised Nash bargaining (FLNS) and social welfare (FSW) models. 

\begin{equation}
    \Omega = \underset{f}{\sum}\pi_{f}\qquad \label{eq:naiveapp}
\end{equation}

\section{Case studies}\label{CaseGT}
For the purpose of this paper two case studies from an industrial liquid market will be examined, a duopoly and an oligopoly comprised of three firms. The examined time horizon is fifteen years discretised into quarterly intervals, $p=1,\cdots,60$. Three different contracts are examined with varying duration, 5 year, 3 year and 1 year contracts $k=1,2,3$ respectively for each of which, we have examined 3 key terms, $s=1,2,3$. The computational experiments were carried in an Intel® Core™i9-10900K CPU @ 3.70GHZ machine using GAMS v41.3 and the MILP solver Gurobi v9.5.2 using 20threads. For the social welfare models the solution time was less than 5 minutes for both case studies. In contrast, for the Nash bargaining models, in the duopoly case study the CPU time was almost an hour, while in the oligopoly case study 10 hours. 
\subsection{Duopoly}
In the duopoly case study the formation of a coalition comprised of two firms $f=2$, $c=97$ customers with $t=315$ tanks and $i=2$ trading products has been evaluated. Out of the 97 customers only 81 are allocated in the status quo market, hence 17 are free allowing for a market growth, and the products are liquid oxygen and liquid nitrogen. The initial market share in terms of profit at the status quo market is 78/22 \% as indicated by Table \ref{tab:Dmark}. From the same Table it can be observed that the social welfare scheme results in a significant profit increase for the smaller firm, i.e. 141\% compared to the status quo, while only 7\% for the case of Firm A. Such a disproportional profit increase for the duopoly firms results in a shift in the market share to 60/40\% respectively. Even though Firm A is still "leading" in the duopoly, the profit share allocation is less favourable for Firm A. In contrast, the Nash bargaining approach results in market share allocation that represents better the stability in existing market dynamics, 69/31 \%. At the same time, Firm A meets a profit expansion of 19\% along with 83\% increase for Firm B. 
\begin{table}
\centering
\caption{Market share and profit increase for different fairness schemes over Status Quo in the duopoly case study.}
\label{tab:Dmark}
 \begin{tabular}{c c c c } 
 \hline
     & \textbf{Status Quo}  & \textbf{Social Welfare} & \textbf{Nash bargaining}\\ \hline
  \textbf{\%Market share A} & 78 & 60    &69 \\ 
   \textbf{\%Market share B} & 22  & 40      & 31\\ 
     \textbf{\%Profit change A} & - & +7    &+19 \\ 
   \textbf{\%Profit change B} & -  & +141      & +83\\ 
 \hline
\end{tabular}
\end{table}

\begin{figure}
    \centering
    \includegraphics[width=1\linewidth]{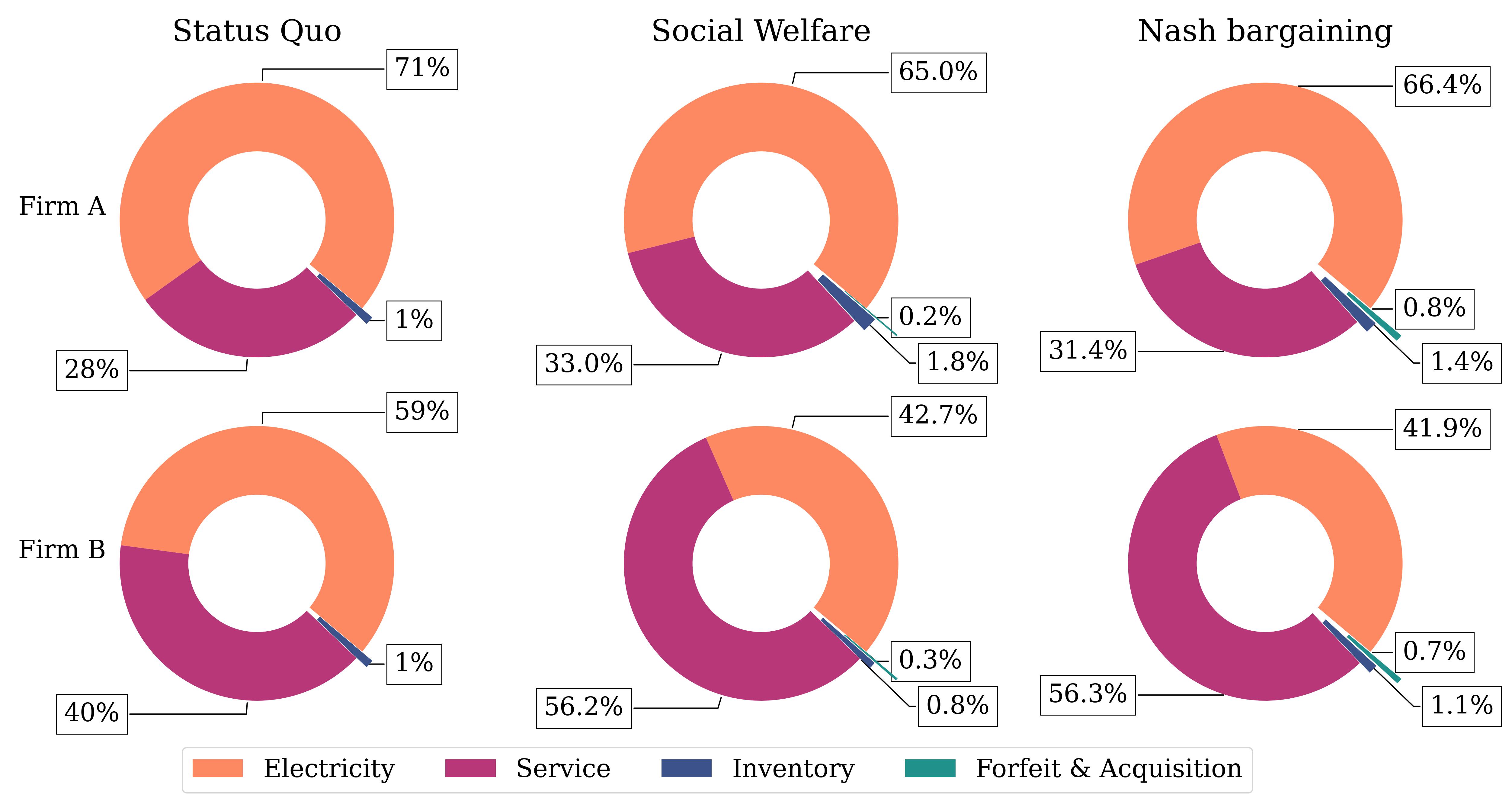}
    \caption{Cost breakdown for duopoly firms and different fairness schemes.}
    \label{fig:Ddonut}
\end{figure}

Apart from the profit allocation, an indicator of the supply chain planning is the cost allocation. Fig.\ref{fig:Ddonut} represents the cost allocation as \% of the the total cost of each firm under the different fairness schemes. For Firm A in the status quo market the major contribution in the cost stems from the Electricity, 71\%, followed by the Service cost, 28\% and only 1\% of inventory. The incorporation of the 17 new customers does not affect the cost breakdown for Firm A as there is a minor variance for social welfare and Nash bargaining schemes, their main difference over the status quo is an additional forfeit and acquisition cost. Examining the second row of Fig.\ref{fig:Ddonut}, it can be observed that the dominance of the electricity cost in the status quo market is replaced by the dominance of the Service cost for both fairness schemes. This shift can be attributed to the market expansion of Firm B. It is noteworthy that for the status quo and social welfare market the demand of the allocated customers is served by product amount by the contracted firm's in-house production, while in the case of the Nash bargaining model swaps between firms cover less than 1\% of the demand. Both firms operate their ASU units  close to the capacity bounds, but even if the ratio of contracted demand to plant capacity is greater than one, the swaps between firms suffice to satisfy the demand and no product amount is purchased from the spot-marke.

\begin{figure}
    \centering
    \includegraphics[width=1\linewidth]{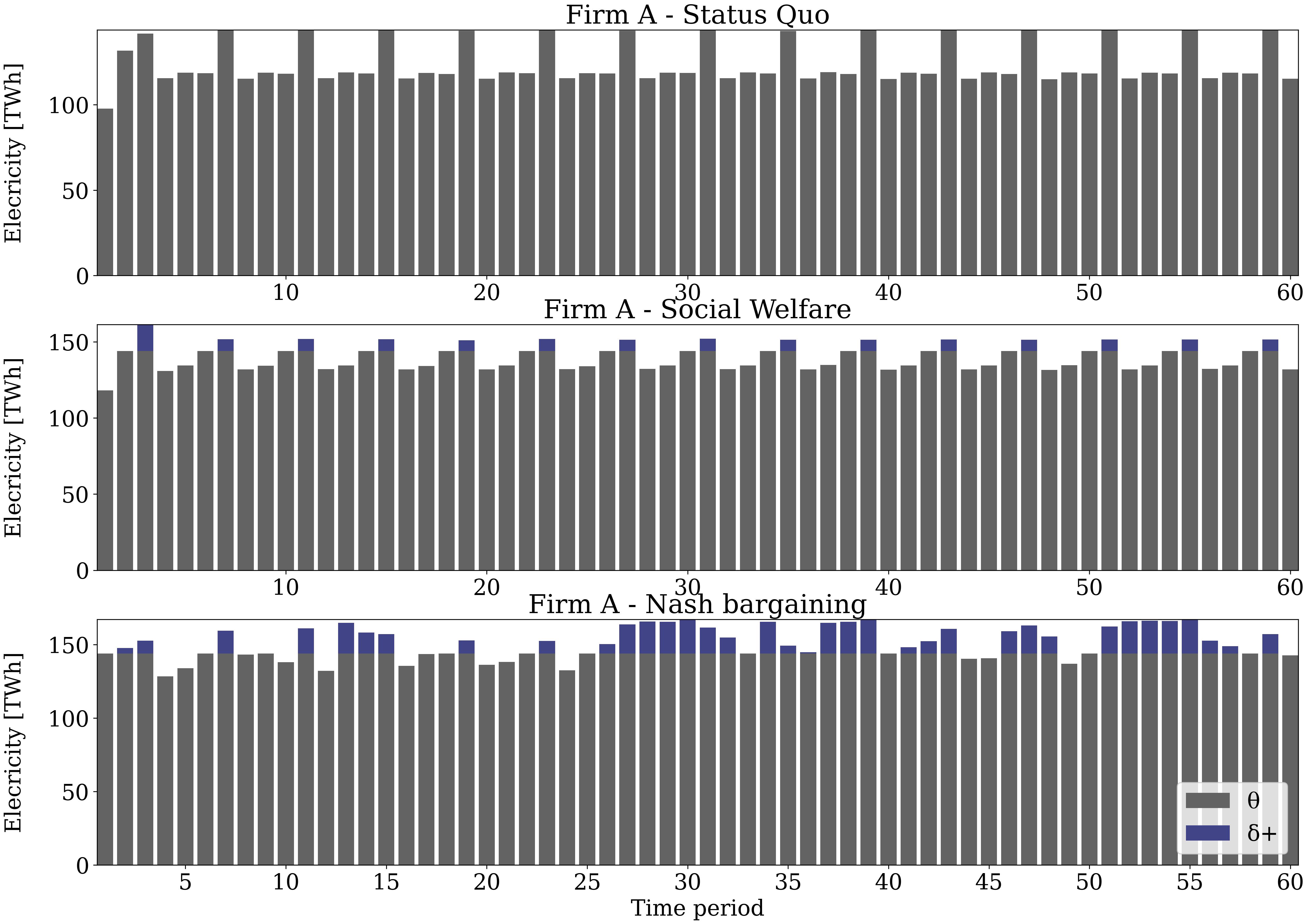}
    \caption{Electricity consumption for different fairness schemes. $\theta$ corresponds to electricity consumed within the pre-agreed threshold with the energy provider while $\delta+$ to electricity used above that threshold.}
    \label{fig:Delec}
\end{figure}

Since the electricity cost is the main cost of Firm A for all of the examined models, Fig.\ref{fig:Delec} focuses on the electricity consumption for the different models. For the status quo case, a periodical variation in the electricity consumption within the pre-agreed threshold with the energy provider can be discerned. In the social welfare scheme a periodicity is also present, however every 4 time periods the electricity consumption exceeds the agreed threshold resulting in higher electricity price. The over-consumption of energy in the Nash bargaining model takes place for the majority of time periods and ranges from 5 to 20 Twh of overpriced electricity. A similar electricity consumption occurs for Firm B. The addition of the 17 customers in the duopoly results in an increased production, hence the re-negotiation of the electricity thresholds with the energy provider will result in reduced electricity cost. The recurrence of the electricity consumption can be attributed to the periodicity of the demand, which in the summer period meets a 20\% increase compared to the other seasons.

\begin{figure}
    \centering
    \includegraphics[width=0.75\linewidth]{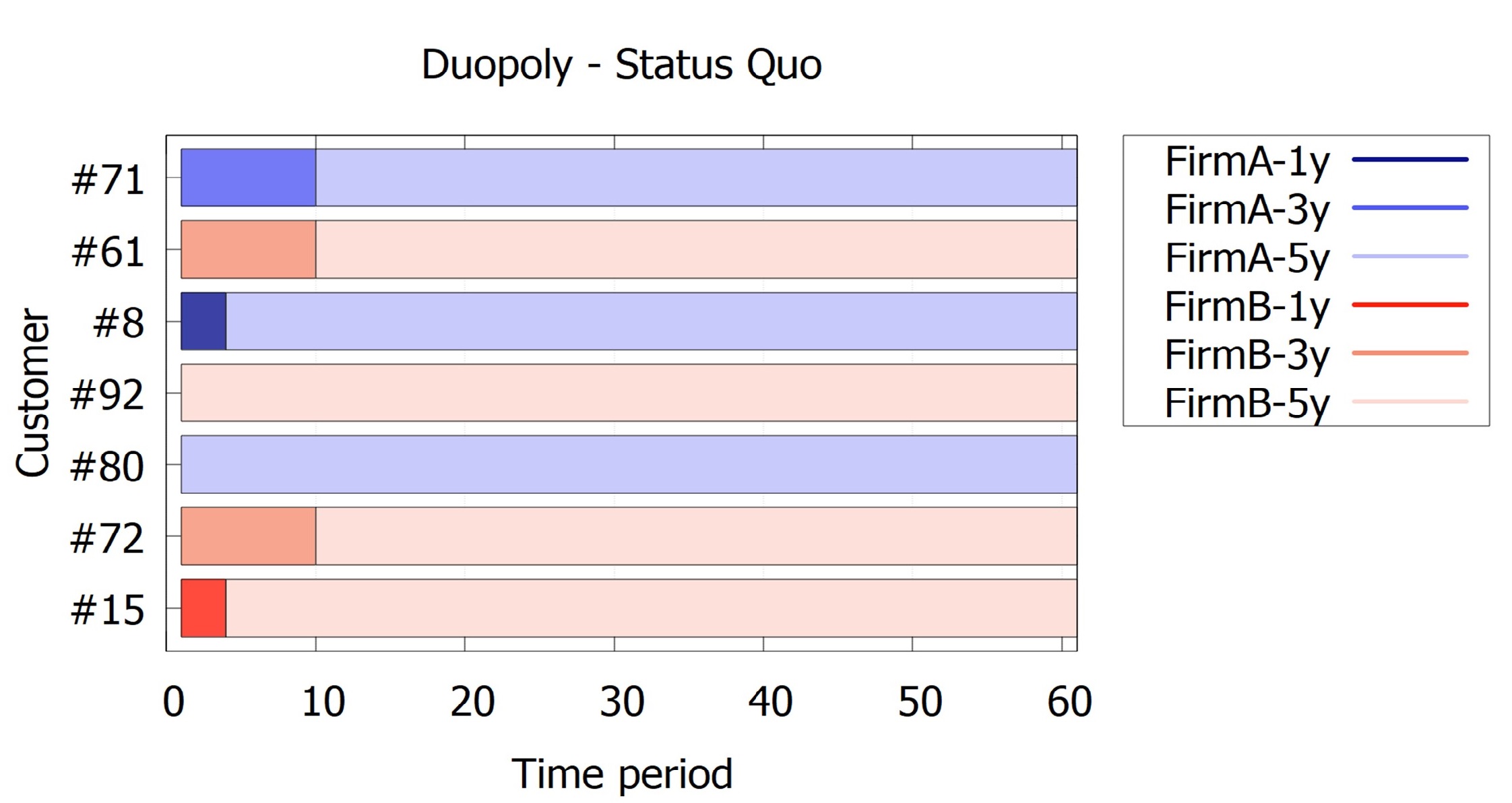}
    \caption{Gantt chart for contract allocation of the customers with the highest demand in the duopoly at the Status Quo market.}
    \label{fig:DSQgantt}
\end{figure}

\begin{figure}
    \centering
    \includegraphics[width=0.75\linewidth]{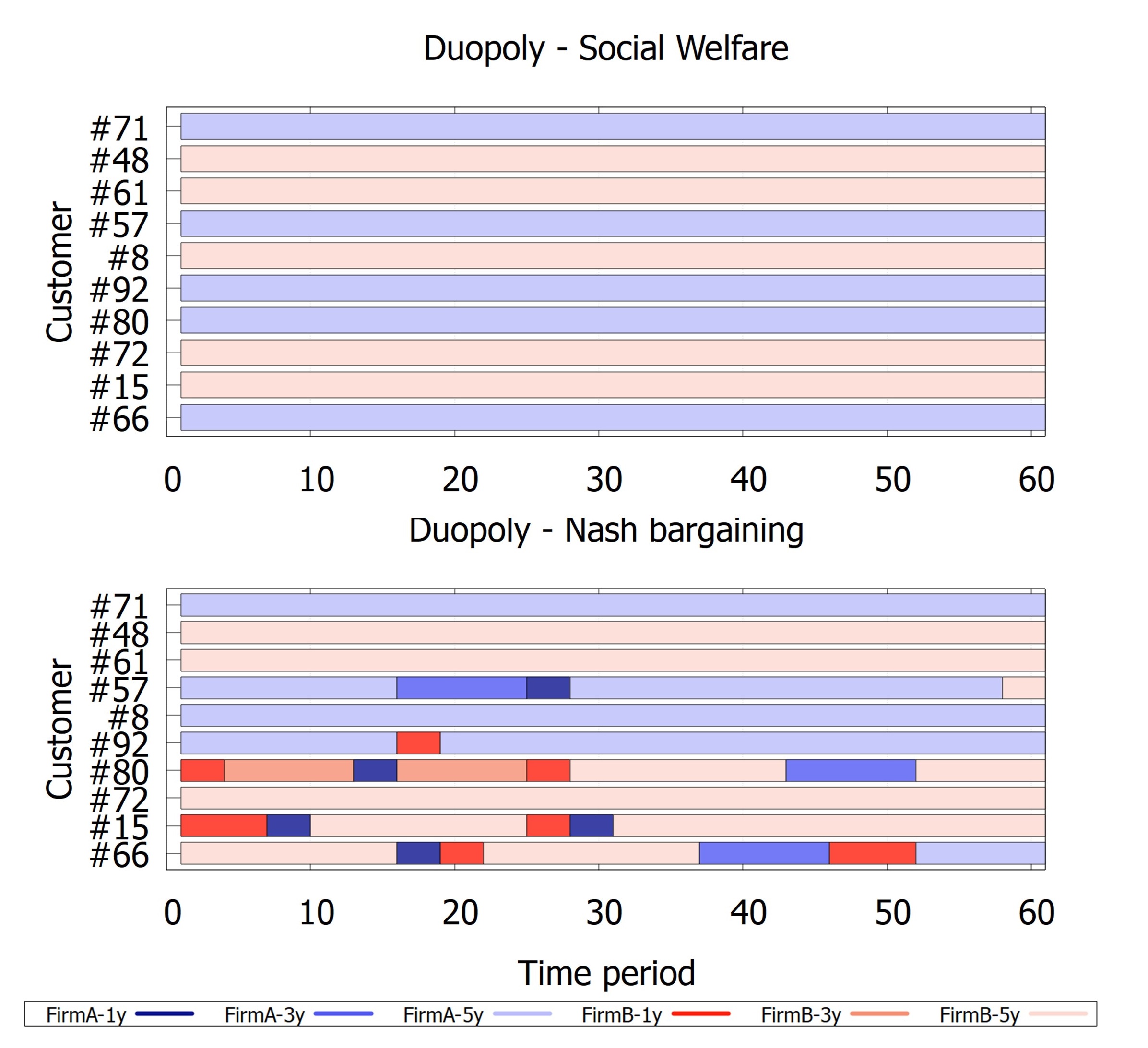}
    \caption{Gantt chart for contract allocation of the customers with the highest demand in the duopoly for the Social Welfare and Nash bargaining models.}
    \label{fig:Dfairgantt}
\end{figure}

The impact of the fairness scheme on the contract and customer allocation can be visualised with the use of the Gantt charts in Figs.\ref{fig:DSQgantt}, \ref{fig:Dfairgantt}. The selected customers for the Gantt charts correspond to those with the highest demand in the market and are displayed in descending order, i.e. customer \#71 has the greatest total demand in the duopoly followed by customer \#48 and so on. At first, it can be noted that customers \#48, \#57 and \#66, even though they belong to the top ten customers of the duopoly, in the status quo market they are not served by any firm. Under a competitive framework, both firms would have aimed to attract these free customers. However, as it can be observed by Fig.\ref{fig:Dfairgantt} the optimal allocation by both fairness schemes is that \#48 is served by Firm B, \#57 by Firm A , while for customer \#66 either Firm A, social welfare model, or a combination of both, Nash bargaining is selected. For the top 3 customers of the duopoly both fairness schemes select the same firm and 5 year contracts. Both in the status quo and social welfare allocation the biggest customers remain with the same firm for the examined time horizon. In contrast, the Nash bargaining scheme results in a higher mobility between firms while employing 1 year and 3 year contracts.The higher customer mobility inflicts a corresponding forfeit and acquisition cost to the firms (see Fig. \ref{fig:Ddonut}). 

Apart from the different duration the contracts have also different base and escalation factors. As a general remark, the base ($\beta_{ictfk}$) decreases with the contract duration, while the escalation factor ($\epsilon_{fpsk})$) decreases, $\beta_{ictf,k=1}<\beta_{ictf,k=2}<\beta_{ictf,k=3}$ and $\epsilon_{fps,k=1}>\epsilon_{fps,k=2}>\epsilon_{fps,k=3}$, where k=1 corresponds to 1 year contract, k=2 to 3 year and k=3 to 5 year contracts respectively.

\subsection{Oligopoly}
 An additional case study has been examined to evaluate the formation of an oligopoly comprised of 3 firms. In this case liquid argon is introduced in the market as well, so the total number of products $i=3$, the total number of customers is $c=81$ and the tanks $t=119$. Out of the 81 customers only 68 are allocated in the status quo market. As it is suggested by Table \ref{tab:Omark} Firm C is dominating in the oligopoly market, holding 47\% of the market share, while Firms A and B have 25\% and 28 \% of the market's profit respectively. In the social welfare scheme the market share between Firm B and C is reversed and Firm B becomes the main supplier of the market. The fact that Firm B becomes the "leader" in the supply chain is achieved by 241\% profit increase compared to the status quo market, whilst Firm C has a 22\% profit decrease. Firm A maintains a similar market share along with a 44\% profit growth. The results of profit allocation in the social welfare scheme favour significantly Firm B, however the realisation of such a market allocation could not be accepted by Firm C. In contrast, Nash bargaining scheme yields a lucrative market allocation for all firms which have a 50\% profit increase each, and at the same time maintains the existing market dynamics. The impact of the fairness scheme is significantly highlighted in the oligopoly case study, since a simplistic approach following the social welfare would result in a market re-structure and significant losses in a firm of the oligopoly, while the Nash bargaining approach results in fair profit allocation based on existing market share.

\begin{table}[h]
\centering
\caption{Market share and profit increase for different fairness schemes over Status Quo in the oligopoly case study.}
\label{tab:Omark}
 \begin{tabular}{c c c c } 
 \hline
     & \textbf{Status Quo}  & \textbf{Social Welfare} & \textbf{Nash bargaining}\\ \hline
  \textbf{\%Market share A} & 25 & 21    &25 \\ 
   \textbf{\%Market share B} & 28  & 57      & 28\\ 
   \textbf{\%Market share C} & 47  & 22 & 47\\ 
     \textbf{\%Profit change A} & - & +44    &+50 \\ 
   \textbf{\%Profit change B} & -  & +241      & +50\\ 
   \textbf{\%Profit change C} & -  & -22 & +50\\ 
 \hline
\end{tabular}
\end{table}

\begin{figure}
    \centering
    \includegraphics[width=1\linewidth]{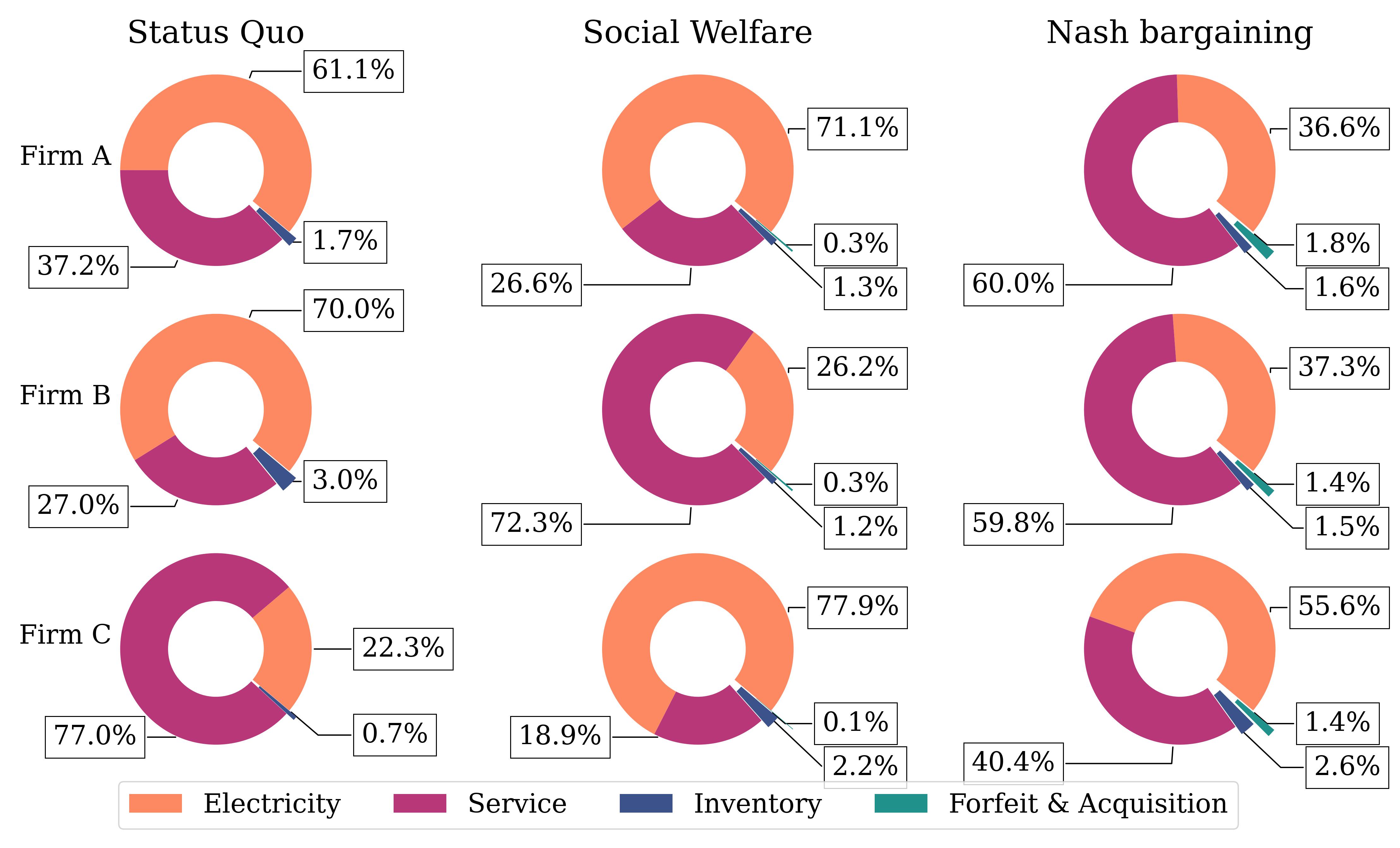}
    \caption{Cost breakdown for oligopoly firms and different fairness schemes.}
    \label{fig:Odonut}
\end{figure}

The differences between the fairness schemes can also be detected on the cost breakdown in Figure \ref{fig:Odonut}. Starting with Firm C, in the status quo market 77\% of the cost is derived from customer service and only 22\% in electricity consumption. In the social welfare model the primal sources of cost are ordered differently, were now the 77\% of  Firm C's total cost stems from electricity consumption and only 19\% from service cost. A more balanced cost allocation occurs in the Nash bargaining model. An analogous alteration is observed for Firm B from a reversed perspective. The cost of Firm B in the status quo market is dominated by the electricity consumption  as 70\%  followed by  27\%  by the service cost, however, this percentages are reversed in the social welfare and Nash bargaining models. Firm A has a similar cost breakdown for the status quo and the social welfare models where 61\% and 71\% of the total cost is attributed to the electricity consumption, here in the Nash bargaining case the Service cost has a greater impact, 60\%, of the total cost. Similar to the duopoly case study, the Nash bargaining approach results in increased forfeit and acquisition cost along with higher inventory costs for all firms, compared to the social welfare model.

\begin{figure}
    \centering
    \includegraphics[width=1\linewidth]{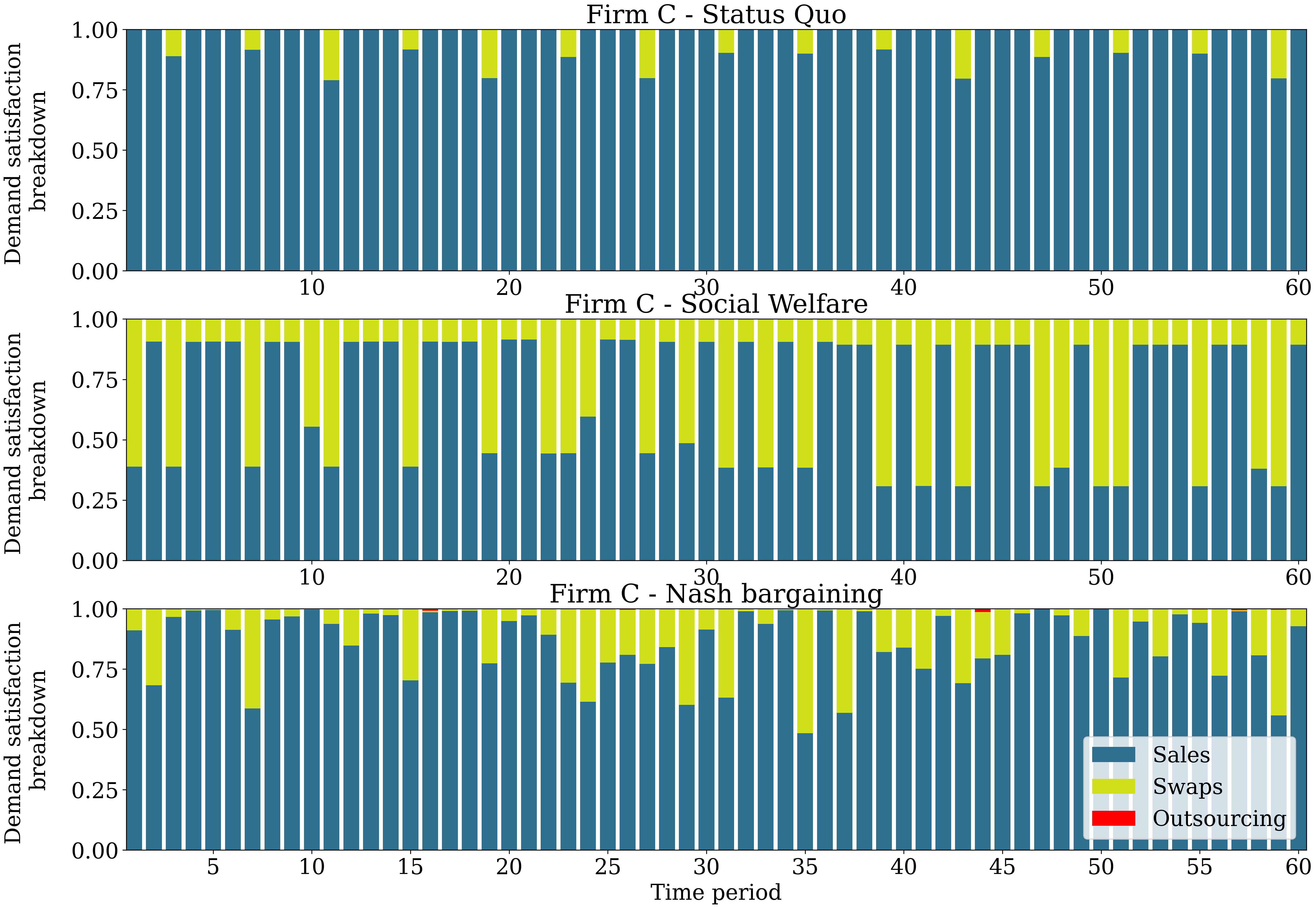}
    \caption{Demand satisfaction breakdown for Firm C and different fairness schemes in the oligopoly case study.}
    \label{fig:Odem}
\end{figure}

To further examine the variance of cost breakdown in Firm C, a more in depth analysis of the service cost is deemed necessary. Figure \ref{fig:Odem} depicts the demand satisfaction breakdown of the customers assigned to Firm C for different fairness schemes. In the status quo market, the demand of Firm C is satisfied mainly by product amount by the contracted firm's in-house production and occasional by swaps from Firm A and B. As Firm C is depreciated to a smaller market share in the social welfare scheme, the dependence to the swaps from the other firms is increased, reaching in half of the time periods a 50\% of product swaps.The operation of the supply chain with Firm C relying on sush a high level in swaps is highly unlikely, since Firm C is the main supplier in the existing market. In the Nash bargaining scheme, Firm C, while maintaining the same market share as the status quo market, covers the customer demand dominantly by product amount by the contracted firm's in-house production. In this model a small fraction of the demand is covered by outsourcing the production in the spot-market, see periods 16, 26, 44, 57 and 59.

\begin{figure}
    \centering
    \includegraphics[width=0.75\linewidth]{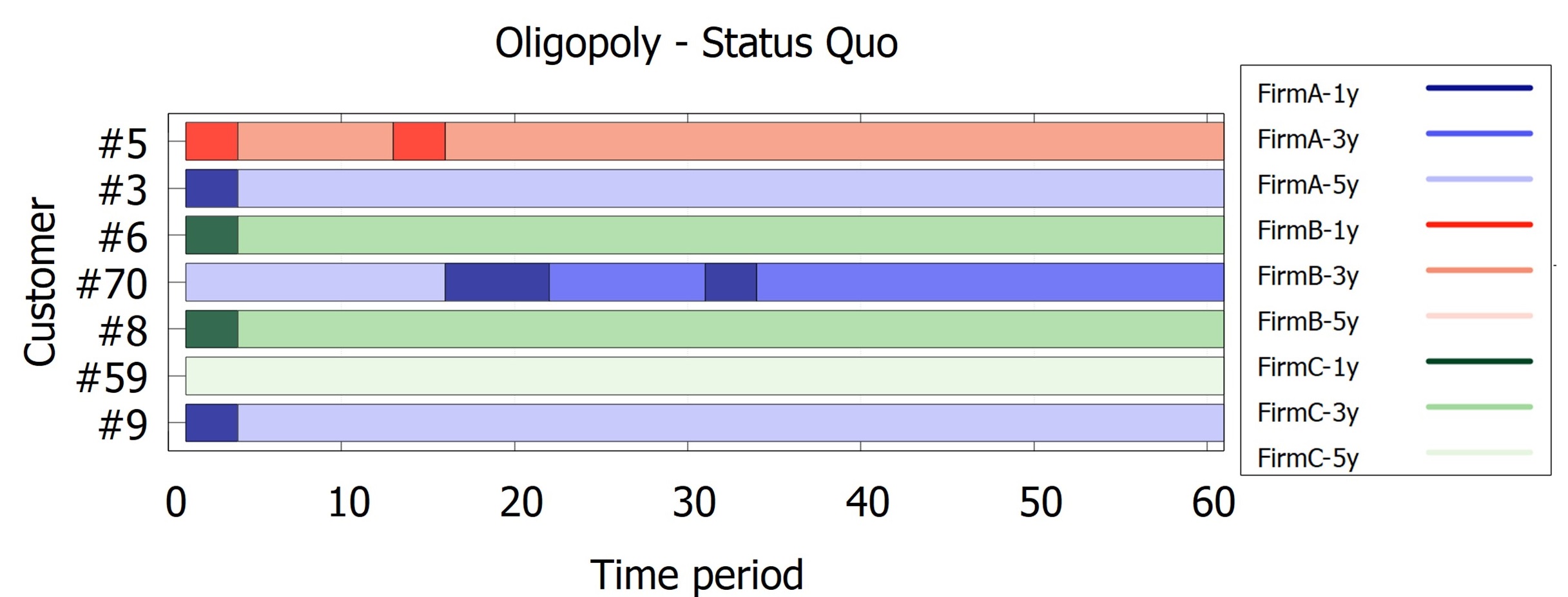}
    \caption{Gantt chart for contract allocation of the customers with the highest demand in the oligopoly at the Status Quo market.}
    \label{fig:OSQgantt}
\end{figure}
\begin{figure}
    \centering
    \includegraphics[width=0.75\linewidth]{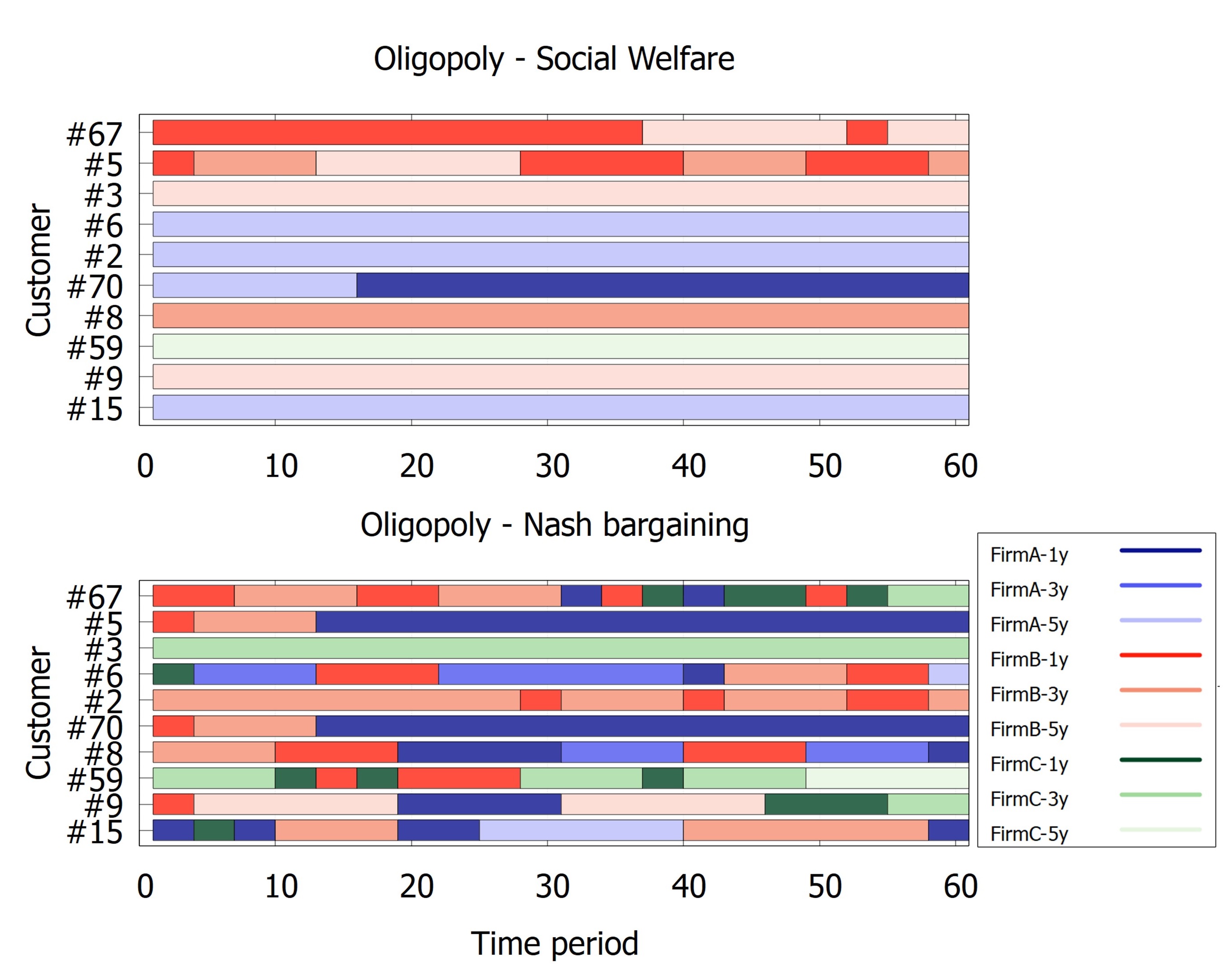}
    \caption{Gantt chart for contract allocation of the customers with the highest demand in the oligopoly for the Social Welfare and Nash bargaining models.}
    \label{fig:Ofairgantt}
\end{figure}

The Gantt charts for the contract and customer allocation are illustrated in Figs. \ref{fig:OSQgantt} and \ref{fig:Ofairgantt}. The examined customers are those with the highest total demand. In this case study, the biggest customer, \#67, is a free customer in the status quo market, so as customer \#15. In the status quo and social welfare market allocation, the customers with the greater demand are fixed in a specific firm for the examined time horizon and the choice of the long term 5 year contracts is prevailing. On the contrary, in the Nash bargaining scheme customers are mainly assigned 1 year or 3 year contracts which facilitates the alteration of customers between firms, e.g. customers \#67, \#6 and \#15 have signed with all three oligopoly firms (Fig.\ref{fig:Ofairgantt}). However, in this case study the exchange of customers between firms is conducted in an excessive and irrational degree at times, looking at customer \#67 from time period 30 the customer changes a supplying firm every year for 4 years.  Despite the increased exchange of customers between firms the corresponding forfeit and acquisition cost corresponds to less than 3\% of the total cost for all firms. To avoid such variation, extra constraints could be introduced preventing such a high customer mobility.

\section{Conclusions}\label{sec:conc}
The problem of fair customer allocation in oligopoly supply chain markets has been addressed in this paper via  two multi-period models. The notion of fairness is expressed via different objective function formulations following two game-theoretic schemes, i.e. social welfare scheme and the Nash bargaining. The former results into an MILP problem while the latter in a non-convex MINLP, which was was approximated with the use of SOS2 variables to a MILP model. The results of the two models were compared to the status quo market which corresponds in an initial customer allocation in which free customers exist, allowing for market growth. The results of the computational experiments suggest that the social welfare  model favours the growth of the smaller firms of the oligopolies at the expense of the bigger firms, customers are signed with contracts of longer duration and remain faithful to a specific firm. In contrast, the Nash bargaining scheme distributes the profit growth among all firms of the oligopolies while maintaining the initial market structure. To achieve that, shorter term contracts are favoured and customers switch firms. 

Future work directions entail the extension of the models so as to take into account the customers as additional players towards a universally fair game. The values of the contract parameters are inherently uncertain, however in this study we have examined a deterministic approach for the contract terms. Taking into account uncertainty in the contract outcomes would provide resilience in the strategic planning of the oligopoly supply chains.

\section*{Acknowledgments}
The authors gratefully acknowledge financial support from EPSRC grants EP/V051008/1, EP/T022930/1 and EP/V050168/1.

\section*{Data and code availability}
The data that has been used for the case studies are confidential. The code of the abstract mathematical model is available upon reasonable request from the corresponding author.

\section*{Appendix}
\begin{figure}[!htb]
\renewcommand{\thefigure}{A1}
    \centering
    \includegraphics[width=0.8\linewidth]{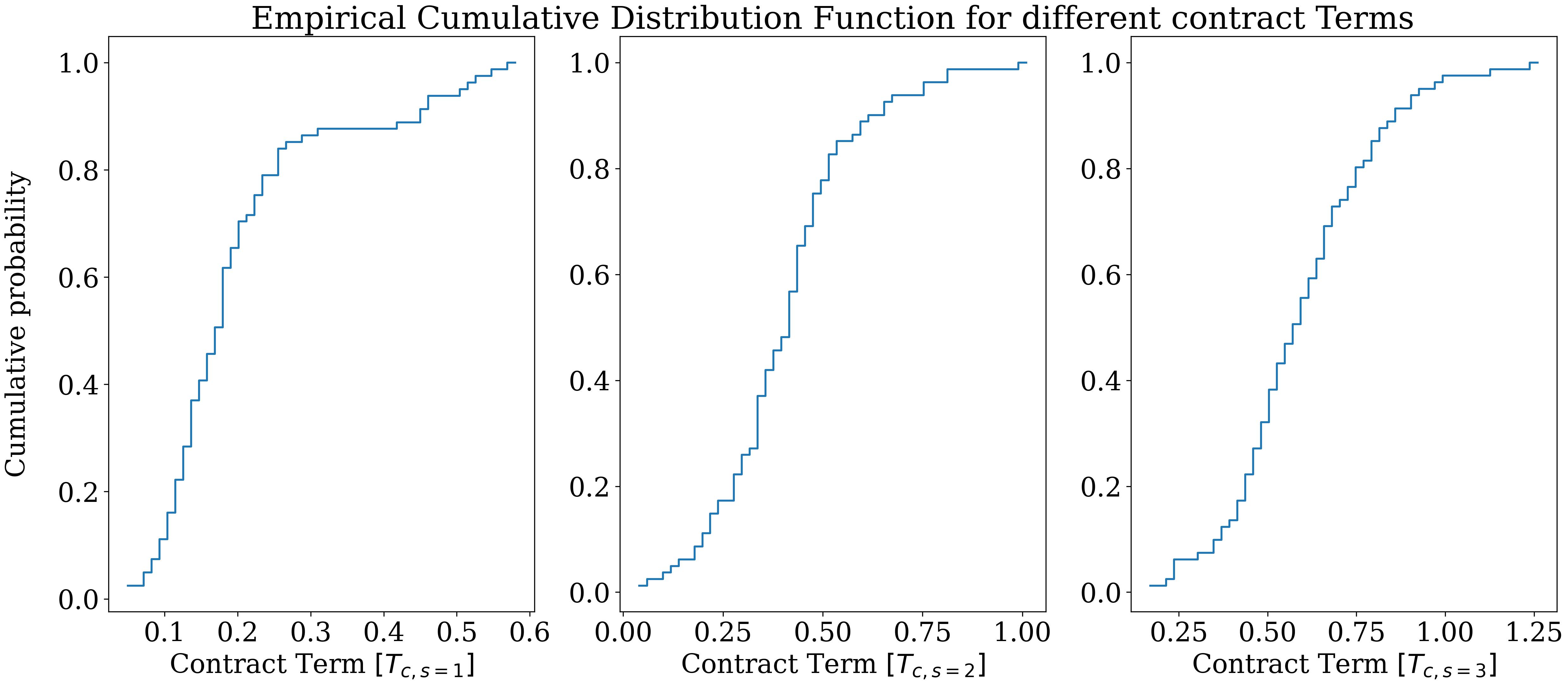}
    \caption{Empirical CDF for for contract terms $T_{cs}$ for the Oligopoly case study}
    \label{fig:Terms}
\end{figure}

\begin{figure}[!htb]
\renewcommand{\thefigure}{A2}
    \centering
    \includegraphics[width=0.8\linewidth]{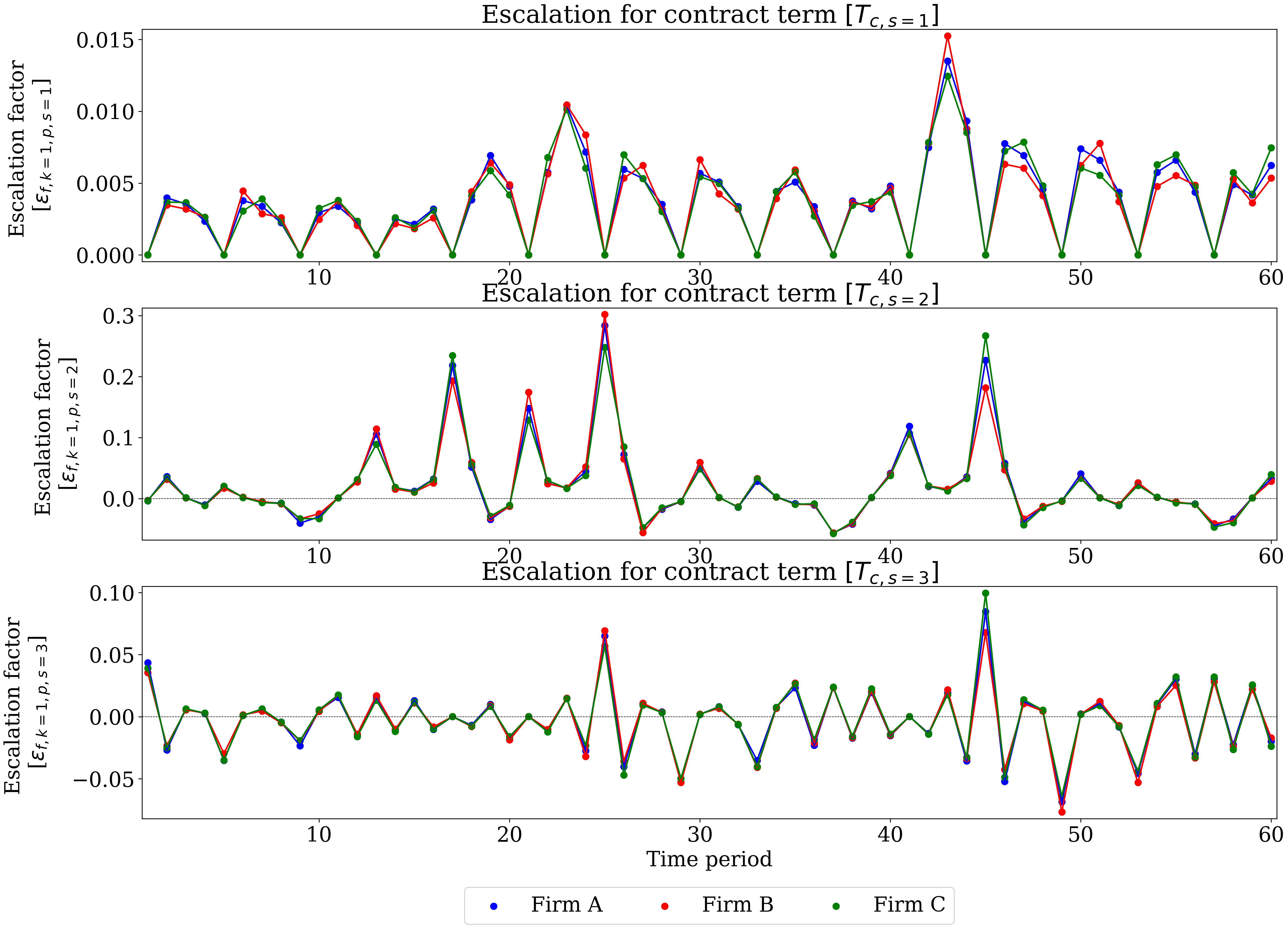}
    \caption{Escalation factors $\epsilon_{fkps}$ for different contract terms and 1 year contracts in the Oligopoly case study}
    \label{fig:escal}
\end{figure}

\bibliographystyle{elsart-harv}

\bibliography{Oligop}

\end{document}